\documentclass[amssymb,amstex,amsmath,reqno,11pt]{amsart}
\usepackage{epsfig}
\input epsf
\usepackage{latexsym}
\usepackage{color}
\usepackage{amsmath}
\usepackage{graphics}
\usepackage{amsthm}
\usepackage{amssymb}

\begin{document}

\title{Stable capillary hypersurfaces  in a wedge}
\author{Jaigyoung Choe and Miyuki Koiso}
\date{}
\thanks{J.C. supported in part by NRF, 2011-0030044, SRC-GAIA;\\\indent M.K. supported in part by
Grant-in-Aid for Scientific Research (B) No. 25287012
of the Japan Society for the Promotion of Science,
and
the Kyushu University Interdisciplinary Programs in Education and Projects in Research Development. }
\address{Jaigyoung Choe: Korea Institute for Advanced Study, Seoul, 130-722, Korea \\
Miyuki Koiso:
Institute of Mathematics for Industry, Kyushu University, 744, Motooka, Nishi-ku, Fukuoka, 819-0395, Japan}
\email{choe@kias.re.kr, koiso@math.kyushu-u.ac.jp}

\maketitle

\newtheorem{theorem}{Theorem}
\newtheorem{cor}{Corollary}
\newtheorem{prop}{Proposition}
\newtheorem{lemma}{Lemma}
\newtheorem{condition}{Condition}
\newtheorem{example}{Example}
\newtheorem{definition}{Definition}
\newtheorem{remark}{Remark}
\newtheorem{conjecture}{Conjecture}
\newtheorem{claim}{Claim}
\newcommand{\rf}[1]{\mbox{(\ref{#1})}}

\renewcommand{\thefootnote}{\fnsymbol{footnote}}



\begin{abstract}
Let $\Sigma$ be a compact immersed stable capillary hypersurface in a wedge bounded by two hyperplanes in $\mathbb R^{n+1}$. Suppose that $\Sigma$ meets those two hyperplanes in constant contact angles and is disjoint from the edge of the wedge. It is proved that if $\partial \Sigma$ is embedded for $n=2$, or if $\partial\Sigma$ is convex for $n\geq3$, then $\Sigma$ is part of the sphere. And the same is true for $\Sigma$ in the half-space of $\mathbb R^{n+1}$ with connected boundary $\partial\Sigma$.
\end{abstract}

\section{Introduction}\label{intro}

The isoperimetric inequality says that among all domains of fixed volume in $\mathbb R^{n+1}$ the one with least boundary area is the round ball. What happens if the boundary area is a critical value instead of the minimum? For this question the more general domains enclosed by the immersed hypersurfaces have to be considered, hence one needs to introduce the oriented volume (as defined in \rf{volume}). Then the answer to the question is that given a compact immersed hypersurface $\Sigma$ in $\mathbb R^{n+1}$, its area is critical among all variations of $\Sigma$ preserving the oriented volume enclosed by $\Sigma$ if and only if $\Sigma$ has constant mean curvature(CMC).

So H. Hopf conjectured that a compact immersed hypersurface of CMC should be a round sphere. To this conjecture W.-Y. Hsiang \cite{H} obtained a counterexample, a CMC immersion of $\mathbb S^3$ in $\mathbb R^4$ which is not round, and in 1986 Wente \cite{W1986} constructed a CMC immersion of a torus in $\mathbb R^3$.

Then, is there an extra condition on a CMC surface $\Sigma$ which guarantees that $\Sigma$ is a sphere? There are some affirmative results in this regard: i) Alexandrov \cite{A} showed that every compact {\it embedded} hypersurface of CMC in $\mathbb R^{n+1}$ is a sphere, ii) Hopf himself \cite{Ho} proved that an immersed CMC {\it 2-sphere} is round, and iii) Barbosa and do Carmo \cite{BD} showed that the only compact immersed {\it stable} CMC hypersurface of $\mathbb R^{n+1}$ is the sphere. A CMC hypersurface $\Sigma$ is said to be stable if the second variation of the $n$-dimensional area of $\Sigma$ is nonnegative for all $(n+1)$-dimensional volume-preserving perturbations of $\Sigma$.

A CMC surface with nonempty boundary along which it makes a constant contact angle with a supporting surface is called a capillary surface. McCuan  \cite{M} and Park \cite{Park2005} proved that an embedded annular capillary surface in a wedge in $\mathbb R^3$ is necessarily part of the sphere. Then the question arises: can one extend the theorems of Alexandrov, Hopf, and Barbosa-do Carmo to the case of capillary surfaces in a wedge or in the half-space? That is, i) show that there is no compact embedded capillary surface of genus $\geq1$ in a wedge (or in the half-space) of $\mathbb R^3$,
ii) is there a compact immersed annular capillary surface of genus 0 (or higher) in a wedge (or in the half-space) which is not part of the sphere? iii) which hypothesis of McCuan's and Park's can be dropped or generalized if the capillary
surface is stable? To question i) McCuan \cite{M}  gave an affirmative answer with the contact angle condition
$\theta_i\leq\pi/2$.
In relation to question ii) Wente \cite{W3} constructed noncompact capillary surfaces bifurcating from the cylinder in a wedge.
 In this paper we have the following answer (\S \ref{main}) to question iii):

\begin{quote}
\em Let $\Sigma$ be a compact immersed stable capillary hypersurface in a wedge bounded by two hyperplanes in $\mathbb R^{n+1}$. Suppose that $\Sigma$ meets those two hyperplanes in constant contact angles and does not hit the edge of the wedge. If $\partial \Sigma$ is embedded for $n=2$, or if $\partial\Sigma$ is convex for $n\geq3$, then $\Sigma$ is part of the sphere. Also, the same conclusion holds if $\Sigma$ is in the half-space of $\mathbb R^{n+1}$ and $\partial\Sigma$ is connected.
\end{quote}

Wente \cite{W1991} simplified Barbosa-do Carmo's proof by using the parallel hypersurfaces and the homothetic contraction. We have found that Wente's method carries over nicely to our capillary hypersurfaces in a wedge and in the half-space. On the other hand, the Minkowski inequality for $\partial\Sigma$ is indispensable in our arguments.

Finally,  it should be mentioned that the stable capillary surfaces in a ball also have been studied very actively. To begin with, Nitsche \cite{N} showed that a capillary disk in a ball $\subset\mathbb R^3$ is a spherical cap. Ros and Souam \cite{RS}
proved that a stable capillary surface of genus 0 in a ball in $\mathbb R^3$ is a spherical cap. They also proved that a stable minimal surface with constant contact angle in a ball $\subset\mathbb R^3$ is a flat disk or a surface of genus 1 with at most three boundary components. Moreover, Ros and Vergasta \cite{RV} showed that a stable minimal hypersurface in a ball $B\subset\mathbb R^n$ which is orthogonal to $\partial B$ is totally geodesic, and that a stable capillary surface in a ball $\subset\mathbb R^3$ and orthogonal to $\partial B$ is a spherical cap or a surface of genus 1 with at most two boundary components.

We would like to thank Professor Monika Ludwig for referring us to the Alexandrov-Fenchel inequality.

\section{Preliminaries}\label{setting}

Let $\Pi_1$ and $\Pi_2$ be two hyperplanes in ${\mathbb R}^{n+1}$ containing the $(n-1)$-plane $\{x_n=0, x_{n+1}=0\}$ and making angles $\alpha,-\alpha$ ($0<\alpha<\pi/2$) with the horizontal hyperplane $\{x_{n+1}=0\}$, respectively.
Let $\Omega\subset\{x_n>0\}$ be the wedge-shaped domain bounded by $\Pi_1$ and $\Pi_2$.
We denote by $\overline{\Omega}$ the closure of $\Omega$.
Denote by $X:(\Sigma , \partial \Sigma )\rightarrow (\overline{\Omega},\partial \Omega)$ an immersion of an
$n$-dimensional oriented compact connected $C^\infty$ manifold $\Sigma$ with nonempty boundary
into  $\Omega$ such that $X(\Sigma^\circ) \subset \Omega$ and
$X(\partial\Sigma)\subset\partial\Omega$, where $\Sigma^\circ :=\Sigma -\partial\Sigma$.
The $(n-1)$-plane $\Pi_0:=\Pi_1\cap\Pi_2=\{x_n=0, x_{n+1}=0\}$
is called the edge of the wedge $\Omega$. In this paper we are concerned only with the immersed surfaces $X(\Sigma)$ which connect $\Pi_1$ to $\Pi_2$ without intersecting $\Pi_0$.

For the immersion $X:(\Sigma , \partial \Sigma )\rightarrow (\overline{\Omega},\partial \Omega)$ the $n$-dimensional area $\mathcal{H}^n(X)$ is written as
\begin{equation*}
\mathcal{H}^n(X)=\int_\Sigma dS
\end{equation*}
where $dS$ is the volume form of $\Sigma$ induced by the immersion $X$. The $(n+1)$-dimensional {\it oriented volume} $V(X)$ enclosed by $X(\Sigma)$ is defined by
\begin{equation}\label{volume}
V(X)=\frac{1}{n+1}\int_\Sigma \langle X,\nu\rangle dS
\end{equation}
where $\nu$ is the Gauss map of $X$.
Here, the Gauss map is defined in the following manner.
The Gauss map $\nu$ is the unit normal vector field along $X$ of which the orientation is determined as follows.
Let $\{e_1, \cdots, e_n\}$ be an oriented frame on
the tangent space $T_p(\Sigma)$, $p \in \Sigma$.
Then $\{dX_p(e_1), \cdots, dX_p(e_n), \nu\}$ is a frame of $\mathbb R^{n+1}$ with positive orientation.

\begin{center}
\includegraphics[width=2.5in]{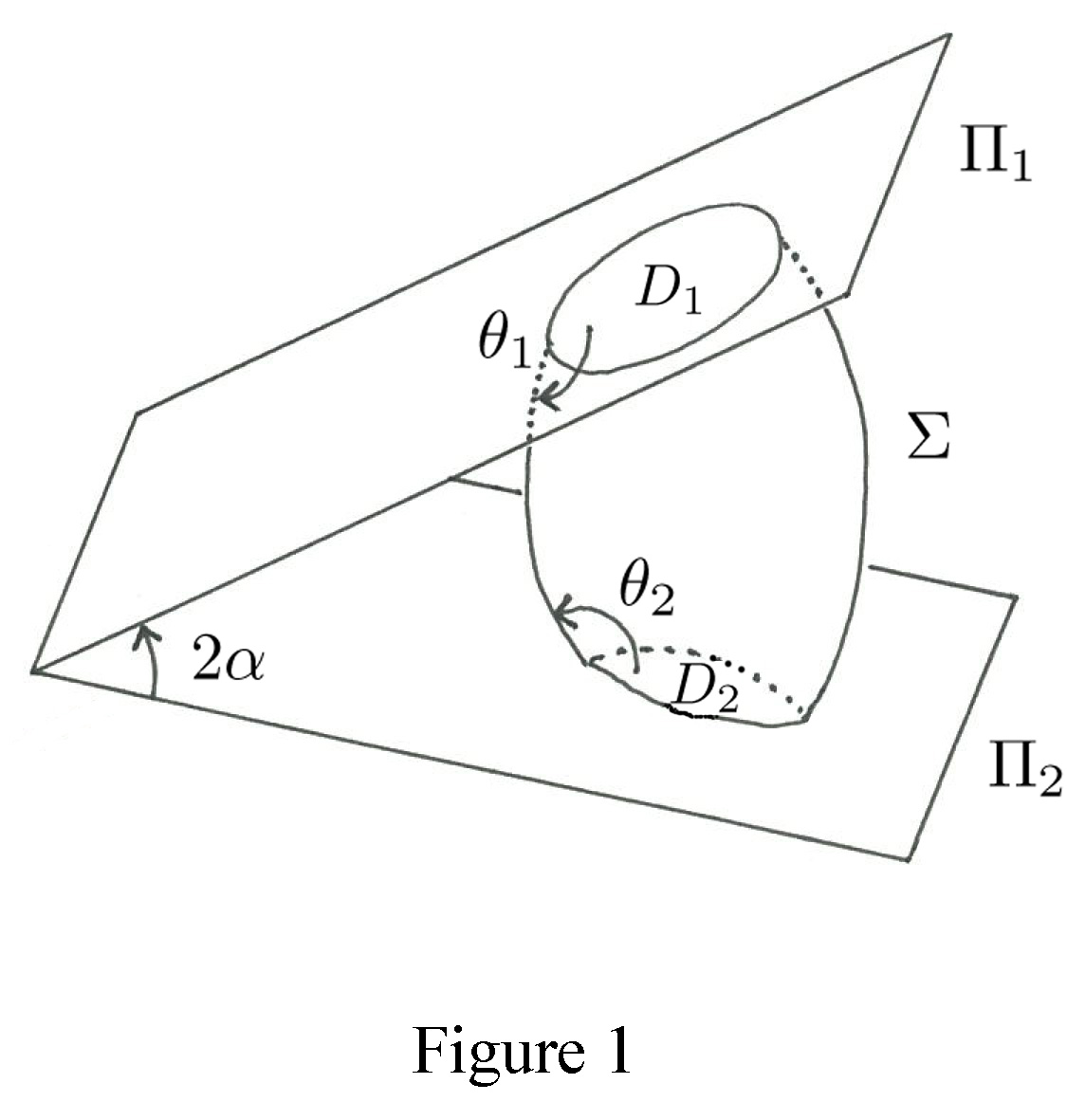}\\
\end{center}

In this paper $X(\Sigma)$ is immersed while $X(\partial\Sigma)$ is assumed to be embedded. And $X(\partial\Sigma)$ will have some influence on the area $\mathcal{H}^n(X)$ through the {\it wetting energy}. Set $C_i=X(\partial\Sigma)\cap\Pi_i$ and let $D_i\subset\Pi_i$ be the domain bounded by $C_i$.
The wetting energy ${\mathcal W}(X)$ of $X$ is defined by
\begin{equation*}\label{w-energy}
{\mathcal W}(X)=\omega_1\mathcal{H}^n(D_1)+\omega_2\mathcal{H}^n(D_2),
\end{equation*}
where $\omega_i$ is a constant with $|\omega_i|<1$ and $\mathcal{H}^n(D_i)$ is the $n$-dimensional area of $D_i$. Then we define the {\it total energy} $E(X)$ of the immersion $X$ by
\begin{equation*}
E(X)=\mathcal{H}^n(X)+\mathcal W(X).
\end{equation*}
Note that $\Sigma\cup D_1\cup D_2$ is a piecewise smooth hypersurface without boundary.
We can extend $\nu:\Sigma \to S^n$ to the Gauss map $\nu:\Sigma \cup D_1 \cup D_2 \to S^n$.
Since the origin of $\mathbb R^{n+1}$ is on the edge $\Pi_0$ of $\Omega$,
$\langle X,\nu\rangle=0$ on $D_1\cup D_2$. Hence
the oriented volume
\begin{equation}\label{widehat}
\widehat{V}(X)=\frac{1}{n+1}\int_{\Sigma\cup D_1\cup D_2}\langle X,\nu\rangle dS
\end{equation}
coincides with $V(X)$.

Let $X_t:(\Sigma , \partial \Sigma )\rightarrow (\overline{\Omega},\partial \Omega)$ be a 1-parameter family of immersions with $X_0=X$. It is well known \cite{KP2007} that  a necessary and sufficient condition for $X$ to be a critical point of the total energy for all variations $X_t$ for which
the volume $\widehat{V}(X_t)$ is constant  is that the immersed surface have constant mean curvature $H$ and that the contact angle $\theta_i$ of $X(\Sigma)$ with $\Pi_i$ (measured between $X(\Sigma)$ and $D_i$) be constant along $C_i$ (see Figure 1). More precisely,
\begin{equation*}\label{contact}
\cos\theta_i=-\omega_i\,\,\,{\rm on}\,\,C_i.
\end{equation*}
The hypersurface $X(\Sigma)$ of constant mean curvature with constant contact angle along $C_i$ will be called a {\it capillary} hypersurface.
A capillary hypersurface is said to be stable if the second variation of $E(X_t)$ at $t=0$ is nonnegative for all volume-preserving perturbations $X_t:(\Sigma,\partial\Sigma)\rightarrow(\overline{\Omega},\partial\Omega)$ of $X(\Sigma)$.


A capillary hypersurface $X(\Sigma)$ in $\overline{\Omega}$ has a nice property called the {\it balancing formula} (\cite{C},\cite{KKS}):

\begin{lemma}
\begin{equation}\label{balance}
nH\mathcal{H}^n(D_i)=-(\sin\theta_i)\mathcal{H}^{n-1}(C_i), \quad i=1, 2.
\end{equation}
\end{lemma}

\noindent
{\it Proof.}
First we remark the following fact:
Let $\hat\Sigma$ be an $m$-dimensional oriented compact connected $C^\infty$ manifold,
and $Y:\hat\Sigma \to {\mathbb R}^{m+1}$ be a continuous map which is a piecewise $C^\infty$ immersion.
Also let $\hat\nu$ be the Gauss map of $Y$.
Then, by using the divergence theorem, we obtain
$$
\int_{\hat\Sigma} \hat\nu \; dS =0.
$$
Now integrate
\begin{equation*}
\Delta_\Sigma X=nH\nu
\end{equation*}
on $\Sigma$ to get
$$
\sum_{i=1}^2\int_{C_i} \eta \;ds = nH\int_\Sigma \nu \;d\Sigma,
$$
where $\eta$ is the outward-pointing unit conormal to $\partial\Sigma$ on $X$.
Then, use the above remark to obtain
\begin{equation}\label{b1}
\sum_{i=1}^2\int_{C_i} \eta \;ds
=-nH \sum_{i=1}^2\int_{D_i} \nu \;dS.
\end{equation}
Denote by $N_i$ the unit normal to $\Pi_i$ that points outward from $\Omega$.
Denote by $n_i$ the inward pointing unit normal to $C_i$ in $\Pi_i$.
Set
\begin{equation}
\epsilon_i := \left \{
\begin{array}{l}
1, \ {\rm if} \ \nu =N_i \ {\rm on} \ D_i,\\
-1, \  {\rm if} \ \nu =-N_i \ {\rm on} \ D_i.
\end{array}
\right.
\end{equation}
Then from (\ref{b1}) we obtain
$$
\sum_{i=1}^2 \int_{C_i} \Bigl\{(\sin\theta_i)\epsilon_iN_i-(\cos\theta_i)n_i \Bigr\} \;ds
+ \sum_{i=1}^2 nH\mathcal{H}^n(D_i)\epsilon_iN_i=0,
$$
that is
$$
\sum_{i=1}^2 (\sin\theta_i)\epsilon_i \mathcal{H}^{n-1}(C_i)N_i
-\sum_{i=1}^2 (\cos\theta_i)\int_{C_i} n_i \;ds
+ \sum_{i=1}^2 nH\mathcal{H}^n(D_i)\epsilon_iN_i=0.
$$
where $\mathcal{H}^{n-1}(C_i)$ is the $(n-1)$-dimensional area of $C_i$.
By using the above remark again, we obtain
$$
\sum_{i=1}^2 \big\{nH\mathcal{H}^n(D_i)+(\sin\theta_i)\mathcal{H}^{n-1}(C_i)\big\}N_i=0.
$$Since $N_1$, $N_2$ are linearly independent, we obtain the formula (\ref{balance}).
\hfill$\Box$

\vskip0.5truecm

Another tool that will be essential in this paper is the formula for the volume of tubes due to H. Weyl \cite{W2}.
Given an immersion $X$ of a compact oriented $n$-manifold $M$ into $\mathbb R^{n+1}$, let $X_t=X+t\nu$ be the one-parameter family of parallel hypersurfaces to $X$. Thanks to the parallelness of $X_t$ one can easily see that $X_t$ has the same unit normal vector field as $X$ and that the area $\mathcal{H}^n(X_t)$ is a polynomial of degree $n$ in $t$. Namely, if $k_1,\ldots,k_n$ are the principal curvatures of $X$, then\\
\begin{eqnarray}
\mathcal{H}^n(X_t)&=&\int_M\prod^n_{i=1}(1-k_it)dS\nonumber\\
&=&a_0+a_1t+a_2t^2+\cdots+a_nt^n,\label{tube}\\
a_0&=&\mathcal{H}^n(X_0),\nonumber\\
a_1&=&-\int_MnHdS,\nonumber\\
a_2&=&\int_M\sum_{i<j}k_ik_j\,dS,\nonumber\\
a_\ell&=&(-1)^\ell\int_M\sum_{i_1<\cdots<i_\ell}k_{i_1}k_{i_2}\cdots k_{i_\ell}dS.\nonumber
\end{eqnarray}
Moreover, the oriented volume $V(X_t)$ satisfies
\begin{equation*}\label{volumearea}
\frac{d}{dt}V(X_t)=\mathcal{H}^n(X_t).
\end{equation*}
Hence
\begin{equation*}
V(X_t)=v_0+v_1t+v_2t^2+\cdots+v_{n+1}t^{n+1},
\end{equation*}
\begin{equation*}
v_1=a_0,\,\,2v_2=a_1,\ldots.
\end{equation*}

\section{Admissible variations}\label{admissible}

Here we assume that our capillary hypersurface $X:(\Sigma,\partial\Sigma)\rightarrow(\overline{\Omega},\partial\Omega)$ has a nonempty boundary component on each $\Pi_i$, $i=1, 2$. But the case when $\Sigma$ is in the half-space and $\partial\Sigma$ is connected can be treated similarly.

To check the stability of $X$
one needs to deal with its volume-preserving variations $X_t:(\Sigma,\partial\Sigma)\rightarrow(\overline{\Omega},\partial\Omega)$. The specific variation that we are going to use arises from the parallel hypersurfaces
\begin{equation*}
X^1_t=X+t\nu.
\end{equation*}
But $X^1_t$ does not satisfy the boundary condition $X^1_t(\partial\Sigma)\subset\partial\Omega$ unless $\theta_i=\pi/2$. To move the boundary to a desired place in $\partial\Omega$, let's apply a translation
\begin{equation*}
X^2_t(p)=p\,+\,ta
\end{equation*}
for some $a\in\mathbb R^{n+1}$.
The vector $a$ is determined in such a way that
\begin{equation*}
X_t^2\circ X_t^1(\partial\Sigma)\subset\partial\Omega.
\end{equation*}
Clearly such a vector uniquely exists as can be seen in Figure 2.

\begin{center}
\includegraphics[width=2.9in]{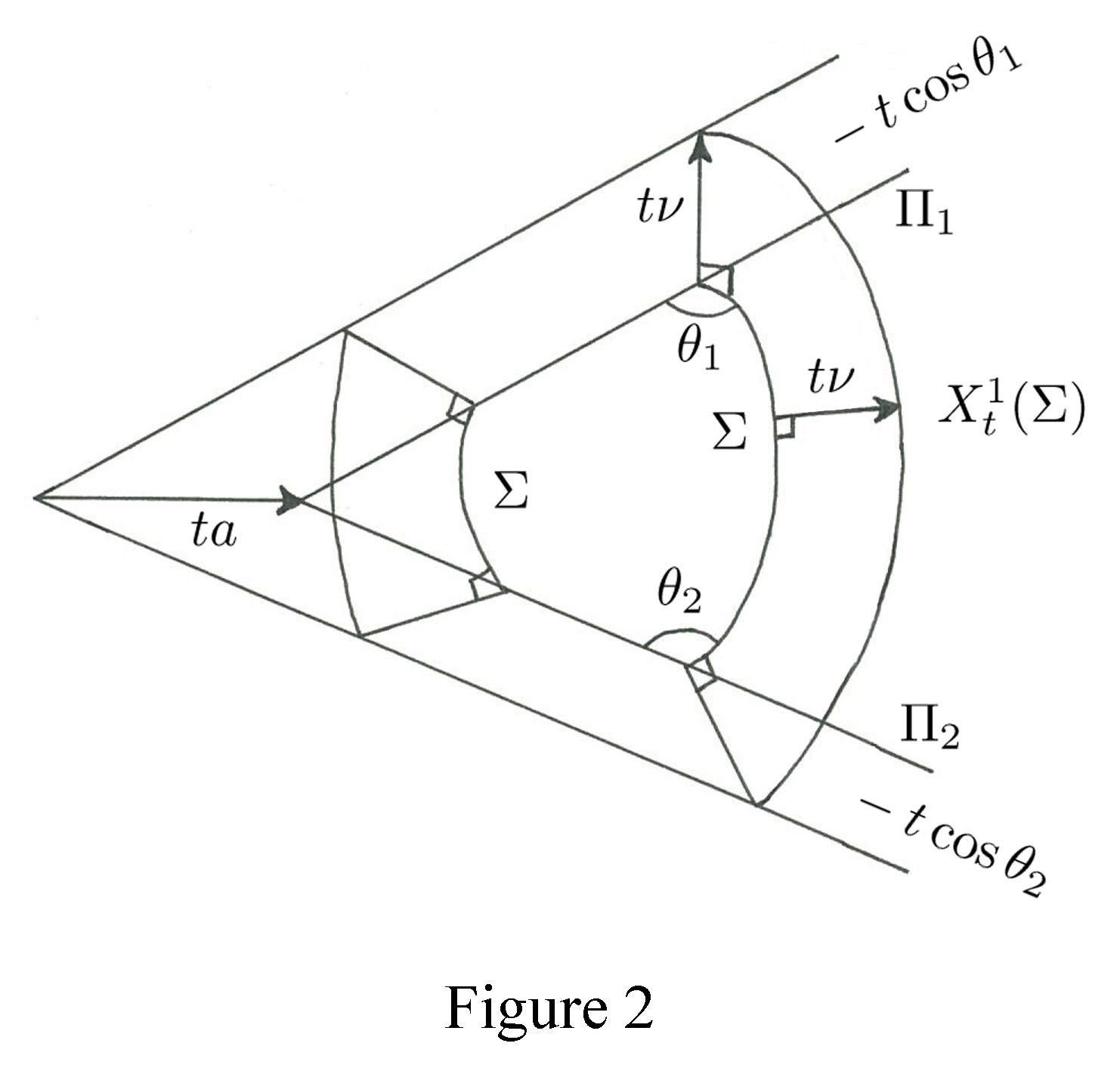}\\
\end{center}

However, $X_t^2\circ X_t^1$ is not volume-preserving. One way of making it into a volume-preserving variation is to deform it by a homothetic contraction:
\begin{equation}\label{st}
X_t:=s(t)X_t^2\circ X_t^1,
\end{equation}
where $s(t)$ satisfies
\begin{equation}\label{v0}
\widehat{V}(X_t)=\widehat{V}(X_0)=v_0.
\end{equation}

In order to compute $\widehat{V}(X_t)$ we first need to consider the oriented volume $\widehat{V}(X_t^2\circ X_t^1)$ enclosed by $X_t^2\circ X_t^1(\Sigma)\cup D_1^t\cup D_2^t $, where $D^t_i\subset\Pi_i$ is the domain bounded by $\Pi_i\cap X_t^2\circ X_t^1(\partial\Sigma)$. Note here that since $X_t^2\circ X_t^1(\Sigma)\cup D_1^t\cup D_2^t $ is closed, the oriented volume $\widehat{V}(X_t^2\circ X_t^1)$ as computed by \rf{widehat} is independent of the translation $X_t^2$. While $t$ increases by $\Delta t$, the oriented volume $\widehat{V}(X_t^2\circ X_t^1)$ increases by $\mathcal{H}^n(X_t^2\circ X_t^1)\Delta t$ on $X_t^2\circ X_t^1(\Sigma)$ and by $-\cos\theta_i\mathcal{H}^n(D_i^t)\Delta t$ on $D_i^t$. Hence
\begin{equation}\label{dve0}
\frac{d}{dt}\widehat{V}(X_t^2\circ X_t^1)=\mathcal{H}^n(X_t^2\circ X_t^1)-\sum_i\cos\theta_i\mathcal{H}^n(D^t_i).
\end{equation}
Calling $-\sum_i\cos\theta_i\mathcal{H}^n(D^t_i)$ the wetting energy $\mathcal{W}(X_t^2\circ X_t^1)$ of $X_t^2\circ X_t^1(\Sigma)$, let's define the total energy
\begin{equation*}
E(X_t^2\circ X_t^1)=\mathcal{H}^n(X_t^2\circ X_t^1)+\mathcal{W}(X_t^2\circ X_t^1).
\end{equation*}

The tube formula \rf{tube} for the capillary hypersurface $\Sigma$ yields
\begin{equation*}
\mathcal{H}^n(X_t^2\circ X_t^1)=a_0+a_1t+a_2t^2+\cdots+a_nt^n,
\end{equation*}
\begin{equation*}
a_0=\mathcal{H}^n(\Sigma),\,\,\,
a_1=-nHa_0,\,\,\,
a_2=\int_\Sigma\sum_{i<j}k_ik_j\,dS,
\end{equation*}
\begin{equation}\label{dve}
\frac{d}{dt}\widehat{V}(X_t^2\circ X_t^1)=E(X_t^2\circ X_t^1).
\end{equation}
Recall $C_i=X(\partial\Sigma)\cap\Pi_i$. Since $X_t^2\circ X_t^1(\Sigma)$ has constant contact angle with $\partial\Omega$ for all $t$, $X_t^2\circ X_t^1(C_i)$ are the parallel hypersurfaces of $p_{\Pi_i}(X^2_t(C_i))$, where $p_{\Pi_i}$ denotes the projection of $\mathbb R^{n+1}$ onto $\Pi_i$. Also recall $\partial D_i=C_i, D_i=D_i^0$. The distance between $X_t^2\circ X_t^1(C_i)$ and $p_{\Pi_i}(X^2_t(C_i))$ is $t\sin\theta_i$. Hence again by the tube formula
for $\mathcal{H}^{n-1}(X_t^2\circ X_t^1(C_i))$, we obtain
\begin{eqnarray*}
\mathcal{H}^n(D_i^t)&=&\mathcal{H}^n(D_i)+\mathcal{H}^{n-1}(C_i)\,t\sin\theta_i-\frac{1}{2}\left(\int_{C_i}(n-1)\bar{H}d\bar{S}\right)t^2\sin^2\theta_i+\\
&&\cdots+ (-1)^{n-1}\frac{1}{n}\left(\int_{C_i}\bar{k}_1\bar{k}_2\cdots\bar{k}_{n-1}d\bar{S}\right)t^{n}\sin^{n}\theta_i,
\end{eqnarray*}
where $\bar{H}$ and $\bar{k}_i$ are the mean curvature and principal curvature of $C_i$ in $\Pi_i$
with respect to the outward unit normal, respectively, and $d\bar{S}$ is the $(n-1)$-dimensional volume form of $C_i$.

Then \rf{dve0} gives
\begin{eqnarray*}
\frac{d}{dt}\widehat{V}(X_t^2\circ X_t^1)&=&a_0-\sum_i\cos\theta_i\mathcal{H}^n(D_i)\nonumber\\
&&-\left(nHa_0+\sum_i\cos\theta_i\sin\theta_i\,\mathcal{H}^{n-1}(C_i)\right)t\nonumber\\
&&+\left(\int_\Sigma\sum_{i<j}k_ik_j\,dS+\frac{1}{2}\sum_i\cos\theta_i\sin^2\theta_i\int_{C_i}(n-1)\bar{H}d\bar{S}\right)t^2\nonumber\\
&&+\cdots.
\end {eqnarray*}
Hence if we write
\begin{equation*}
E(X_t^2\circ X_t^1)=e_0+e_1t+\cdots+e_nt^n,
\end{equation*}
\rf{dve} yields
\begin{eqnarray}
e_0&=&a_0-\sum_i\cos\theta_i\mathcal{H}^n(D_i),\nonumber\\
e_1&=&-nHa_0-\sum_i\cos\theta_i\sin\theta_i\,\mathcal{H}^{n-1}(C_i),\label{eee}\\
e_2&=&\int_\Sigma\sum_{i<j}k_ik_j\,dS+\frac{1}{2}\sum_i\cos\theta_i\sin^2\theta_i\int_{C_i}(n-1)\bar{H}d\bar{S}.\nonumber
\end{eqnarray}

On the other hand, if we let
\begin{equation*}
\widehat{V}(X_t^2\circ X_t^1)=v_0+v_1t+v_2t^2+\cdots+v_{n+1}t^{n+1},
\end{equation*}
then it follows from \rf{st}, \rf{v0} and the binomial series that
\begin{eqnarray*}
s(t)^n&=&v_0^\frac{n}{n+1}(v_0+v_1t+v_2t^2+\cdots+v_{n+1}t^{n+1})^{-\frac{n}{n+1}}\\
&=&1-\frac{n}{n+1}\left(\frac{v_1}{v_0}\right)t+\left\{\frac{n(2n+1)}{2(n+1)^2}
\left(\frac{v_1}{v_0}\right)^2-\frac{n}{n+1}\left(\frac{v_2}{v_0}\right)\right\}t^2+\cdots.
\end{eqnarray*}
Thus
\begin{eqnarray}
E(X_t)&=& s(t)^nE(X_t^2\circ X_t^1(\Sigma))\nonumber\\
&=&e_0+\left\{e_1-\frac{n}{n+1}\left(\frac{v_1}{v_0}\right)e_0\right\}t\label{eprime}\\
&&+\left\{e_2-\frac{n}{n+1}\left(\frac{v_1}{v_0}\right)e_1+\frac{n(2n+1)}{2(n+1)^2}\left(\frac{v_1}{v_0}
\right)^2e_0-\frac{n}{n+1}\left(\frac{v_2}{v_0}\right)e_0\right\}t^2\nonumber\\
&&+\cdots.\nonumber
\end{eqnarray}
From \rf{dve} we have
\begin{equation}\label{v1}
v_1=e_0,\,\,2v_2=e_1,
\end{equation}
and the fact that $E^\prime(0)=0$ in \rf{eprime} implies
\begin{equation}\label{v01}
v_0=\frac{n}{n+1}\,\frac{e_0^2}{e_1}.
\end{equation}
Substituting the identities of \rf{v1} and \rf{v01} into the coefficient of $t^2$  in \rf{eprime} yields
\begin{equation*}
E^{\prime\prime}(0)/2=\frac{1}{2ne_0}\left\{2ne_0e_2-(n-1)e_1^2\right\}.
\end{equation*}
Hence from \rf{eee} we get
\begin{eqnarray*}
ne_0E^{\prime\prime}(0)&=&2n\left(a_0-\sum_i\cos\theta_i\mathcal{H}^n(D_i)\right)\times\\
&&\left(\int_\Sigma\sum_{i<j}
k_ik_jdS + \frac{1}{2}\sum_i\cos\theta_i\sin^2\theta_i\int_{C_i}(n-1)\bar{H}d\bar{S}\right)\\
&&-(n-1)\left(nHa_0 + \sum_i\cos\theta_i\sin\theta_i\,\mathcal{H}^{n-1}(C_i)\right)^2.
\end{eqnarray*}
Then the balancing formula \rf{balance} yields
\begin{eqnarray*}
\left(nHa_0 + \sum_i\cos\theta_i\sin\theta_i\,\mathcal{H}^{n-1}(C_i)\right)^2=n^2H^2\left(a_0-\sum_i\cos\theta_i\mathcal{H}^n(D_i)\right)^2.
\end{eqnarray*}
Therefore
\begin{eqnarray}
ne_0E^{\prime\prime}(0)&=&\left(a_0-\sum_i\cos\theta_i\mathcal{H}^n(D_i)\right)\times\nonumber\\
&&\left(2n\int_\Sigma\sum_{i<j}k_ik_jdS+n\sum_i\cos\theta_i\sin^2\theta_i\int_{C_i}(n-1)\bar{H}d\bar{S}\right.\nonumber\\
&&\left.-\int_\Sigma n^2(n-1)H^2dS+n^2(n-1)H^2\sum_i\cos\theta_i\mathcal{H}^n(D_i)\right)\nonumber\\
&=&\left(a_0-\sum_i\cos\theta_i\mathcal{H}^n(D_i)\right)\times\nonumber\\
&&\left(-\int_\Sigma\sum_{i<j}(k_i-k_j)^2dS+n\sum_i\cos\theta_i\sin^2\theta_i\int_{C_i}(n-1)\bar{H}d\bar{S}\right.\nonumber\\
&&\left.+n^2(n-1)H^2\sum_i\cos\theta_i\mathcal{H}^n(D_i)\right)\nonumber\\
&=&\left(a_0-\sum_i\cos\theta_i\mathcal{H}^n(D_i)\right)\times\nonumber\\
&&\left\{-\int_\Sigma\sum_{i<j}(k_i-k_j)^2dS\right.\label{second}\\
&&\left.+(n-1)\sum_i\cos\theta_i\sin^2\theta_i\left(n\int_{C_i}
\bar{H}d\bar{S}+\frac{\mathcal{H}^{n-1}(C_i)^2}{\mathcal{H}^n(D_i)}\right)\right\},\nonumber
\end{eqnarray}
where the balancing formula \rf{balance} is used again in the last equality.

We shall see in the next section that
$$n\int_{\partial D_i}
\bar{H}d\bar{S}+\frac{\mathcal{H}^{n-1}(\partial D_i)^2}{\mathcal{H}^n(D_i)}\geq0.$$

\section{Theorem}\label{main}
We are now ready to state the  theorem of this paper.\\

\begin{theorem}
Let $W$ be a wedge in $\mathbb R^{n+1}$ bounded by two hyperplanes $\Pi_1,\Pi_2$. And let $\Sigma\subset W$ be a compact oriented immersed hypersurface that is disjoint from the edge $\Pi_1\cap\Pi_2$ of $W$, having embedded boundary $\partial\Sigma\subset\Pi_1\cup\Pi_2$, and satisfying $\partial\Sigma\cap\Pi_i=\partial D_i$ for a nonempty bounded domain $ D_i$ in $\Pi_i$. Suppose that $\Sigma$ is a stable capillary hypersurface in $W$. In other words, $\Sigma$ is an immersed constant mean curvature hypersurface making a constant contact angle $\theta_i\geq\pi/2$ with $D_i$ such that for all volume-preserving perturbations {\rm (}for the oriented volume enclosed by $\Sigma\,\cup D_1\cup D_2${\rm )} the second variation of the total energy
 $$E(\Sigma)=\mathcal{H}^{n}(\Sigma)-\cos\theta_1 \mathcal{H}^{n}(D_1)-\cos\theta_2 \mathcal{H}^{n}(D_2)$$
  is nonnegative. \\
1) If $n=2$, then $\Sigma $ is part of the 2-sphere.\\
2) If $n\geq3$ and $D_1,D_2$ are convex, then $\Sigma$ is part of the $n$-sphere.

Conversely, if $\Sigma$ is part of the $n$-sphere, then it is stable.

Moreover, the same conclusion holds when $\Sigma$ is in the half-space of $\mathbb R^{n+1}$ and $\partial\Sigma$ is connected.
\end{theorem}

\noindent {\it Proof.}
We prove the theorem for $\Sigma$ in a wedge, and the proof for $\Sigma$ in the half-space is similar.

When $n=2$, \rf{second} becomes
\begin{eqnarray*}
&&2e_0E^{\prime\prime}(0)=\left(a_0-\sum_i\cos\theta_i\mathcal{H}^2(D_i)\right)\times\\
&&\left\{-\int_\Sigma(k_1-k_2)^2dS+\sum_i\cos\theta_i\sin^2\theta_i\left(2\int_{\partial D_i}
kds+\frac{\mathcal{H}^{1}(\partial D_i)^2}{\mathcal{H}^2(D_i)}\right)\right\},
\end{eqnarray*}
where $k$ is the geodesic curvature of $\partial D_i$
with respect to the outward unit normal along $\partial D_i$. Note that on the smooth Jordan curve $\partial D_i$, $\int_{\partial D_i}kds=-2\pi$. Hence the isoperimetric inequality of $D_i$ and the angle condition $\cos\theta_i\leq0$ yield
$$E{''}(0)\leq0.$$
Therefore $\Sigma$ needs to be umbilic everywhere if it is stable.

When $n\geq3$, Minkowski (p.1191, \cite{O}) showed that for a convex domain $D\subset\mathbb R^n$ with mean curvature $H$ on $\partial D$,
\begin{equation*}\label{minkowski}
n\int_{\partial D} |H|dS\leq\frac{\mathcal{H}^{n-1}(\partial D)^2}{\mathcal{H}^n(D)}.
\end{equation*}
Hence it follows from \rf{second}  that the stable $\Sigma$ is all umbilic.

If $\Sigma$ is part of the $n$-sphere, then $\Sigma$ is the minimizer of the energy $E$ among all embedded hypersurfaces in $\Omega$ enclosing the same volume (\cite{ZAT}).
The proof is similar to that of Theorem 4.1 in \cite{KP2007}; the method is essentially the same as in \cite{W1967}.
Hence $\Sigma$ is stable for all $n\geq2$.
\hfill$\Box$

\begin{remark} {\rm
Our contact angle condition $\theta_i\geq\pi/2$ is quite natural because McCuan \cite{M} proved the nonexistence of embedded capillary surfaces with $\theta_i\leq\pi/2$ in a wedge of $\mathbb R^3$. Also it had been experimentally observed that a wedge forces the liquid drops(bridges) with $\theta_i\leq\pi/2$ toward its edge.
}\end{remark}

\begin{remark} {\rm
Theorem holds not only in a wedge $W$ but also in
a domain bounded by $n+1$ hyperplanes of which the normals are linearly independent and all of which pass a common fixed point,
as pointed out by Sung-Ho Park.
 }\end{remark}

\section{Minkowski's inequality}\label{minko}
The Minkowski inequality is not well known among geometers and its proof is not easily available in the literature. So in this section let's give a sketchy proof of the Minkowski inequality. First we need to introduce the mixed volume \cite{S}.

The {\it Minkowski sum} of two sets $A$ and $B$ in $\mathbb R^n$ is the set
$$A+B=\{a+b\in\mathbb R^n:a\in A,b\in B\}.$$
Given convex bodies $K_1,\ldots,K_r$ in $\mathbb R^n$, the volume of the Minkowski sum $\lambda_1K_1+\cdots+\lambda_rK_r,\lambda_i\geq0,$ of the scaled convex bodies $K_i$ is a homogeneous polynomial of degree $n$ given by
$$\mathcal{H}^n(\lambda_1K_1+\cdots+\lambda_rK_r)=\sum_{j_1,\ldots,j_n=1}^rV(K_{j_1},\ldots,K_{j_n})
\lambda_{j_1}\cdots\lambda_{j_n}.$$
$V(K_{j_1},\ldots,K_{j_n})$  is called the {\it mixed volume} of $K_{j_1},\ldots,K_{j_n}$. The mixed volume is uniquely determined by the following three propeties:\\ \indent i) $V(K,\ldots,K)=\mathcal{H}^n(K)$,\,\,\, ii) $V$ is symmetric,\,\,\, iii) $V$ is multilinear.\\ A remarkable property of the mixed volume is the {\it Alexandrov-Fenchel inequality}:
$$V(K_1,K_2,K_3,\ldots,K_n)^2\geq V(K_1,K_1,K_3,\ldots,K_n)\cdot V(K_2,K_2,K_3,\ldots,K_n).$$

For a convex body $K\subset\mathbb R^n$ and a unit ball $B\subset\mathbb R^n$, the mixed volume
$$W_j(K):=V(\overbrace{K,K,\ldots,K}^{n-j\,{\rm times}},\overbrace{B,B,\ldots,B}^{j\,{\rm times}})$$
is called the $j$-th {\it quermassintegral} of $K$.
The Steiner formula says that the quermassintegrals of $K$ determine the volume of the parallel bodies of $K$:
$$\mathcal{H}^n(K+tB)=\sum_{j=0}^n\left(\begin{array}{c}n\\j\end{array}\right)W_j(K)t^j.$$

Comparing the Steiner formula for a convex domain $D\subset\mathbb R^n$ with its tube formula, one can obtain
\begin{eqnarray*}
W_0(D)&=& \mathcal{H}^n(D),\\
nW_1(D)&=& \mathcal{H}^{n-1}(\partial D),\\
nW_2(D)&=&\int_{\partial D}|H|dS,\\
n(n-1)(n-2)W_3(D)&=& 2\int_{\partial D}\sum_{i<j}k_ik_jdS.\\
\end{eqnarray*}
The Alexandrov-Fenchel inequality for the quermassintegrals yields
\begin{eqnarray*}
W_1(D)^2&\geq& W_0(D)W_2(D),\\
W_2(D)^2&\geq& W_1(D)W_3(D).
\end{eqnarray*}
Consequently,
\begin{eqnarray}
n\int_{\partial D}|H|dS&\leq& \frac{\mathcal{H}^{n-1}(\partial D)^2}{\mathcal{H}^n(D)},\label{Mink}\\
\int_{\partial D}\sum_{i<j}k_ik_jdS&\leq& \frac{(n-1)(n-2)}{2}\,\frac{\left(\int_{\partial D}|H|dS\right)^2}{\mathcal{H}^{n-1}(\partial  D)}\nonumber\\
&\leq&\frac{(n-1)(n-2)}{2n^2}\,\frac{\mathcal{H}^{n-1}(\partial  D)^3}{\mathcal{H}^n(D)^2},\label{16}
\end{eqnarray}
where \rf{Mink} is the desired Minkowski inequality. \\

\begin{remark}{\rm
\rf{Mink} is the isoperimetric inequality when $D$ is a domain in $\mathbb R^2$ and so is \rf{16} when $D\subset\mathbb R^3$, because $$\int_{\partial D\subset\mathbb R^2}|k|ds=2\pi\,\,\,\,{\rm and}\,\,\,\int_{\partial D\subset\mathbb R^3}k_1k_2\;dS=4\pi.$$
}\end{remark}

\begin{remark}{\rm
Let $D_t\subset\mathbb R^n$ be the parallel domain with distance $t$ to $D$. Then \rf{Mink} is equivalent to
\begin{equation*}
n\,\frac{\mathcal{H}^{n-1}(\partial  D_t)'}{\mathcal{H}^{n-1}(\partial D_t)}\leq\frac{(n-1)\mathcal{H}^n(D_t)'}{\mathcal{H}^n(D_t)},
\end{equation*}
or equivalently,
\begin{equation*}
\left(\frac{\mathcal{H}^{n-1}(\partial  D_t)^n}{\mathcal{H}^n(D_t)^{n-1}}\right)'\leq0.
\end{equation*}
Hence  the isoperimetric quotient $\mathcal{H}^{n-1}(\partial  D_t)^n/\mathcal{H}^n(D_t)^{n-1}$ decreases as $t$ increases. Indeed, the parallel domain $D_t$ becomes rounder and rounder as $t$ increases.
}\end{remark}


\end{document}